\documentclass[10pt]{amsart}
\usepackage{amsmath}
\usepackage{amssymb}
\usepackage{amsthm}

\usepackage[all]{xy}

\newtheorem{theorem}{Theorem}[section]
\newtheorem{lemma}[theorem]{Lemma}
\newtheorem{corollary}[theorem]{Corollary}

\newtheorem{proposition}[theorem]{Proposition}

\newtheorem{example}[theorem]{Example}
\theoremstyle{definition}
\newtheorem{remark}[theorem]{Remark}
\newtheorem{definition}[theorem]{Definition}

\newenvironment{equationth}{\stepcounter{theorem}\begin{equation}}{\end{equation}}

\def\Z{\mathbb Z}
\def\C{\mathbb C}

\def\0{\underline 0}

\begin{document}

\title{The Cayley--Bedregal reduction for Milnor and Hirzebruch--Milnor classes}

\author{R. Callejas-Bedregal}
\address{Roberto Callejas-Bedregal\\ Centro de Ci\^encias Exatas e da Natureza, Universidade Federal da Para\'{i}ba-UFPb, Jo\~ao Pessoa, PB - Brasil.}
\thanks{Roberto Callejas-Bedregal passed away before this paper was completed.}

\author{F. Caraballo}
\address{Franquiz Caraballo\\ Instituto de Matem\'aticas-Cuernavaca, Universidad Nacional Aut\'onoma de M\'exico.}
\email{franquiz.caraballo@im.unam.mx}

\author{M. F. Z. Morgado}
\address{Michelle F. Z. Morgado\\ Instituto de Bioci\^encias Letras e  Ci\^encias Exatas.\\Universidade Estadual Paulista-UNESP, Brasil.}
\email{m.morgado@unesp.br}

\author{J. Seade}
\address{Jos\'e Seade\\ Instituto de Matem\'aticas-Cuernavaca, Universidad Nacional Aut\'onoma de M\'exico.}
\email{jseade@im.unam.mx}

\thanks{Research partially supported by CAPES, CNPq and FAPESP 
(Brazil), by PAPIIT-UNAM IN101424 (Mexico) and a grant from Simons Foundation International [SFI-MPS-T-Institutes-00011977-JS]}

\dedicatory{\begin{minipage}{\textwidth}\raggedleft\vspace*{2ex}
Dedicated to the memory of Roberto Callejas-Bedregal,\\
whose ideas were central to this work.
\end{minipage}}



\begin{abstract}
The Cayley--Bedregal reduction replaces a complete intersection
$\{f_1=\dots=f_r=0\}$ by the hypersurface $\sum_i\varphi_if_i=0$ in the
projectivized dual bundle. For an isolated complete intersection
singularity of corank $c$, we express the Milnor class of this
hypersurface on its singular fibre $\mathbb P^{c-1}$ in terms of the
Milnor numbers of generic systems obtained by taking linear
combinations of the defining equations. Its degree is the original
Milnor number,
and the higher components recover the Milnor numbers of the systems
with fewer equations. These numbers are invariants of the analytic germ. 

Our second main result concerns the Hodge-theoretic refinement. For a
complete intersection defined by a regular section of a rank-$r$
algebraic vector bundle, the Hirzebruch--Milnor class is determined by
that of the Cayley--Bedregal hypersurface: the push-forward of the
normalized Hirzebruch--Milnor class of the hypersurface, in the sense
of Cappell, Maxim, Sch\"urmann and Shaneson, is
$(-y)^{r-1}$ times that of the complete intersection, with no
positivity assumption on the bundle. The factor $(-y)^{r-1}$ comes
from the motivic decomposition of the projective bundle, and at $y=-1$
the identity recovers Aluffi's comparison of Milnor classes. When the
tautological bundle is very ample, the hypersurface formula of Maxim,
Saito and Sch\"urmann then gives the Hirzebruch--Milnor class of the
complete intersection through the vanishing cycles of a single
function, avoiding their induction through intermediate complete
intersections and their restriction to powers of one line bundle. In integral homology we prove the analytic form of
Aluffi's identity, and obtain formulas for the Milnor class in terms
of local Milnor fibres and for the Chern--Schwartz--MacPherson class
of any projective scheme.
\end{abstract}
\maketitle

\section*{Introduction}

Milnor classes extend the Milnor number to singular loci of arbitrary
dimension. For a local complete intersection $Z$, they measure the
difference between its virtual Chern class and its
Chern--Schwartz--MacPherson class \cite{Sch,MacP,Ful}. This difference
is supported on the singular locus. When the singularities are
isolated, its signed degree is the sum of the local Milnor numbers
\cite{SS,Suwa}; for a positive-dimensional singular locus, nonzero components may
occur in every dimension up to $\dim\operatorname{Sing}Z$
\cite{BLSS2,BSS}.

The hypersurface case admits formulas in terms of characteristic
cycles, Segre classes, vanishing cycles and L\^e cycles
\cite{Aluffi,Par-Pr,Massey0,Massey,BMS,BMS-err,Max}. Complete
intersections introduce several defining equations and a normal
bundle of higher rank \cite{BMS-3,MSS}. The classical Cayley trick replaces a system of equations by one
hypersurface in a larger space. The \emph{Cayley--Bedregal reduction}
turns this idea into a method for characteristic classes: it reduces the
Milnor and Hirzebruch--Milnor classes of a complete intersection to
those of a hypersurface, and it recovers the local class over an
isolated complete intersection singularity.

Let $M$ be a compact complex manifold of dimension $n$, and let $s$
be a regular section of a holomorphic vector bundle $E\to M$ of rank
$r$. Write $Z=Z(s)$, with inclusion $i:Z\hookrightarrow M$. On the
projective bundle $p:\mathbb P(E^\vee)\to M$, evaluation against $s$
defines a section $\hat s$ of $\mathcal O_{\mathbb P(E^\vee)}(1)$:
\[
\hat s(x,[\varphi])=\varphi(s(x)).
\]
Its zero scheme $X=Z(\hat s)$ is a hypersurface, locally defined by
the Cayley equation $\sum_i\varphi_if_i=0$ when $s=(f_1,\dots,f_r)$.
We introduced this reduction in \cite{CBMS}, following a proposal
of Roberto Callejas-Bedregal: use the Cayley hypersurface to recover
the Milnor classes of the complete intersection. The name
Cayley--Bedregal reduction honours his memory.

Aluffi subsequently adopted this reduction in his formula for the
Chern--Schwartz--MacPherson class of an arbitrary embeddable scheme
\cite[\S3.1]{AluffiCSM}, explicitly crediting \cite{CBMS} for the
construction. For a regular section, he proved that
\[
p_*\mathcal M(X)=i_*\mathcal M(Z)
\]
in the Chow group of $M$
\cite[Corollary~3.6]{AluffiCSM}. We give the analytic form of this
identity in integral homology (Theorem~\ref{t:final}). The argument
follows his: the Chern--Schwartz--MacPherson comparison comes from
the fibres of $X\to M$, and the virtual comparison is a universal
identity of Chern classes. A proof by evaluation on universal
quotient bundles also applies to Hirzebruch classes.

Aluffi's hypersurface formula expresses the Milnor class of $X$
in terms of the Segre class of its singularity subscheme
\cite[Theorem~I.4]{Aluffi}. Our local result identifies the components
of the class on a singular fibre in terms of Milnor numbers of
generic systems of defining equations. Let $x$ be an isolated singular
point of $Z$,
with $c=r-\operatorname{rank}(ds_x)$. The singular locus of $X$ over
$x$ is a projective space $S_x\simeq\mathbb P^{c-1}$. For a generic
flag of systems consisting of $r-k$ linear combinations of the
defining equations, write $\mu_k$ for the corresponding Milnor number,
$0\leq k\leq c-1$. Thus $\mu_0=\mu(Z,x)$. With
$H=c_1(\mathcal O_{S_x}(1))$, the Parusi\'nski--Pragacz class on $S_x$
is
\[
\mathfrak M_x=\mu_0[pt]+
\sum_{k=1}^{c-1}(-1)^k\mu_k(1+H)^{k-1}\cap[\mathbb P^k]
\]
(Theorem~\ref{c:local}). Its zero-dimensional degree is
Parusi\'nski's generalized Milnor number \cite{Par}, giving
$\mu(X,S_x)=\mu(Z,x)$. The higher components recover the other
$\mu_k$, which are invariants of the analytic germ $(Z,x)$
(Remark~\ref{r:invariance}). The degree is computed on a projective
model obtained by finite determinacy. Non-characteristic hyperplane
restriction then determines the higher components. A local smoothing
calculation gives a topological interpretation of the degree. In corank one, the
hypersurface equation is the sum of an equation for $Z$ and a
nondegenerate quadratic form in $2(r-1)$ variables
(Proposition~\ref{p:corank-one}).

Our second main result extends the comparison to Hodge-theoretic classes. Suppose that $M$ is
a nonsingular complete complex algebraic variety and $E$ is algebraic.
Let $T_{y*}$ be the normalized motivic Hirzebruch class
\cite{BSY}, and put
$\mathcal M_{y*}=T^{\rm vir}_{y*}-T_{y*}$, following
\cite{CMSS,MSS}. Then
\[
p_*\mathcal M_{y*}(X)=(-y)^{r-1}i_*\mathcal M_{y*}(Z)
\qquad\text{in }H_{2*}(M,\mathbb Q)[y]
\]
(Theorem~\ref{t:hirzebruch}). This holds for arbitrary $E$, without
an ampleness assumption. The factor $(-y)^{r-1}$ has a geometric
origin. The fibres of $X\to M$ are $\mathbb P^{r-1}$ over $Z$ and
$\mathbb P^{r-2}$ elsewhere, so in the relative Grothendieck group
\[
[X\to M]=[\mathbb P^{r-2}]\cdot[\operatorname{id}_M]
+\mathbb L^{r-1}\cdot[Z\to M],
\qquad\mathbb L=[\mathbb A^1].
\]
The motivic comparison follows from $\chi_y(\mathbb L)=-y$; the
virtual comparison follows by universality. At $y=-1$ their
difference recovers Aluffi's identity with rational coefficients.

Several formulas follow directly from these comparisons. Combining
Aluffi's identity with Parusi\'nski--Pragacz
\cite{Par-Pr} expresses the Milnor class of $Z$ through the Euler
characteristics of the Milnor fibres of the local Cayley equations
(Theorem~\ref{t:closed}). This applies to arbitrary singularities
and arbitrary holomorphic bundles. Dropping regularity gives a
formula for the Chern--Schwartz--MacPherson class of an arbitrary
zero scheme, hence of every projective scheme
(Theorem~\ref{t:arbitrary}). When $M$ is projective, the corrected
L\^e-class formula of \cite{BMS-err} gives a corresponding expression
through the global L\^e classes of $X$ (Corollary~\ref{c:Le}).

For Hirzebruch--Milnor classes, the hypersurface formula of Maxim,
Saito and Sch\"urmann \cite{MSS} yields a formula in terms of the
vanishing cycles of the single function $\hat s/\hat s'$ on the
complement of $Z(\hat s')$, where $s'$ is general. This application
requires $\mathcal O_{\mathbb P(E^\vee)}(1)$ to be very ample
(Theorem~\ref{t:closed-y}). Passing to $X$ avoids the induction
through intermediate complete intersections used in \cite{MSS}
for sections of powers of a fixed very ample line bundle.

Dirks, Maxim and Olano \cite{DMO} define Hirzebruch--Milnor type
classes, including spectral refinements, by deformation to the normal
cone and Verdier--Saito specialization. Their Remark~1.2 asks how
these classes relate to those of the affine Cayley hypersurface
$\sum_i y_if_i=0$ on $M\times(\mathbb A^r\setminus\{0\})$ when the
$f_i$ are global functions. Our comparison concerns normalized
Hirzebruch--Milnor classes after proper push-forward from the
projectivized hypersurface. Relating it to their specialization
construction at the level of mixed Hodge modules and spectral classes
is a further question.

The comparison formula in \cite[Theorem~6.4]{CBMS} must be corrected.
Subsection~\ref{ss:old} gives a counterexample and explains why the
stratified transversality hypothesis introduced in that 2018 version
excludes the singular case for $r\geq2$. Aluffi's identity in \cite {AluffiCSM} supplies
the correct comparison for the construction proposed there.
We thank the referee of that earlier version for pointing out the
problem with that transversality hypothesis.

\section{The comparison identity}

Let $M$ be a compact complex manifold of dimension $n$, let $E\to M$ be a holomorphic vector bundle of rank $r\geq1$, and let $s$ be a regular holomorphic section. Thus its local components form a regular sequence along the zero scheme $Z=Z(s)$, a local complete intersection of dimension $n-r$. Write $i:Z\hookrightarrow M$ for the inclusion. Regularity is assumed throughout except in Subsection~\ref{ss:arbitrary}. 

\subsection{The Cayley--Bedregal reduction}

The projective bundle $p:\mathbb P(E^\vee)\to M$ has tautological sequence
\begin{equation*}
    0 \to \mathcal{O}_{\mathbb{P}(E^\vee)}(-1) \to p^*E^\vee \to Q \to 0
\end{equation*}
which dualizes to
\begin{equation}\label{eq:dual_tautological_seq}
0\rightarrow Q^\vee \rightarrow p^{*}E\rightarrow
{\mathcal O}_{\mathbb{P}(E^{\vee})}(1)\rightarrow 0.
\end{equation}
The image of $p^*s$ is a section $\hat s$ of $\mathcal O_{\mathbb P(E^\vee)}(1)$, given by $\hat s(x,[\varphi])=\varphi(s(x))$. 

\begin{lemma}\label{l:fibers}
Let $q:Z(\hat s)\to M$ be the restriction of $p$, and let $x\in M$.
\begin{enumerate}
\item If $x\in Z$, then $q^{-1}(x)=p^{-1}(x)\simeq\mathbb P^{r-1}$. Otherwise $q^{-1}(x)$ is a hyperplane in $p^{-1}(x)$.
\item The section $\hat s$ is regular, and $Z(\hat s)$ is a hypersurface.
\item In a trivialization $s=(f_1,\dots,f_r)$ over $U\subset M$, the singular points of $Z(\hat s)$ are precisely the pairs $(x,[\varphi])$ with $x\in Z\cap U$ and
\begin{equation*}
    \sum_{i=1}^r\varphi_i\,df_i|_x = 0.
\end{equation*}
In particular, $q({\rm Sing}\,Z(\hat s))={\rm Sing}\,Z$, and
$Z(\hat s)$ is nonsingular if and only if $Z$ is nonsingular.
\end{enumerate}
\end{lemma}

A scheme-theoretic description of the singularity subscheme of
$Z(\hat s)$ is given in \cite[Lemma~3.1]{AluffiCSM}.

\begin{proof}
The equation $\varphi(s(x))=0$ gives the fibre description.
In a local trivialization $s=(f_1,\dots,f_r)$, the induced section is
$\sum_i\varphi_if_i$. It does not vanish identically on any component
of the smooth ambient space, so it is regular.

In the chart $\varphi_r=1$, write
\[
g=f_r+\sum_{i<r}\varphi_if_i,
\qquad
dg=df_r+\sum_{i<r}(\varphi_i\,df_i+f_i\,d\varphi_i).
\]
On $g=0$, the conditions $\partial g/\partial\varphi_i=f_i=0$ for
$i<r$ also imply $f_r=0$. Hence $g=dg=0$ is equivalent to
$x\in Z$ and $\sum_i\varphi_i\,df_i|_x=0$.
Finally, regularity of $s$ implies that $x\in Z$ is singular exactly
when $ds_x$ has rank less than $r$, or equivalently when such a
nonzero $\varphi$ exists.
\end{proof}

\subsection{Chern classes and the comparison}

\begin{definition}
    The Fulton--Johnson and Chern--Schwartz--MacPherson classes of $Z$ are
    \[
    c^{FJ}(Z)=c(TM|_Z)c(E|_Z)^{-1}\cap[Z],\qquad
    c^{SM}(Z)=c_*(\mathbf 1_Z).
    \]
    Here $[Z]$ is the fundamental class with scheme multiplicities,
    $\mathbf 1_Z$ is the characteristic function of its support, and
    $c_*$ is MacPherson's natural transformation, normalized by
    $c_*(\mathbf 1_Y)=c(TY)\cap[Y]$ for nonsingular $Y$ \cite{MacP}.
    We use the Milnor-class convention
    \[
    \mathcal M(Z)=(-1)^{n-r}\bigl(c^{FJ}(Z)-c^{SM}(Z)\bigr).
    \]
    More generally, if a regular section of $G\to W$ defines a local
    complete intersection $Y$ of dimension $d$, put
    $c^{FJ}(Y)=c(TW|_Y)c(G|_Y)^{-1}\cap[Y]$ and
    $\mathcal M(Y)=(-1)^d(c^{FJ}(Y)-c^{SM}(Y))$.
    This applies in particular to $Z(\hat s)$.
\end{definition}
These classes lie in the homology of $Z$; their ambient images will
be indicated by push-forward. In the algebraic setting,
\eqref{f:SM}, \eqref{f:FJ} and \eqref{the formula} are
\cite[Propositions~3.4 and~3.3, Corollary~3.6]{AluffiCSM}.
We follow Aluffi's argument in the analytic setting, using the
functorial Chern--Schwartz--MacPherson transformation in homology.

\begin{lemma}{\cite[Proposition~3.4]{AluffiCSM}}\label{l:SM}
    In $H_{2*}(M,\mathbb Z)$,
    \begin{equation}\label{f:SM}
        i_*c^{SM}(Z) = p_*c^{SM}(Z(\hat s))-(r-1)\,c(TM)\cap[M].
    \end{equation}
\end{lemma}
\begin{proof}
    The push-forward of a constructible function integrates it over
    the fibres with respect to Euler characteristic. Lemma~\ref{l:fibers}(1)
    therefore gives
    \[
    p_*\mathbf 1_{Z(\hat s)}=(r-1)\mathbf 1_M+\mathbf 1_Z.
    \]
    Apply $c_*$, which commutes with proper push-forward, and use
    $c_*(\mathbf 1_M)=c(TM)\cap[M]$.
\end{proof}

\begin{remark}
    Taking degrees in \eqref{f:SM} gives 
    \begin{equation*}
        \chi(Z(\hat s))=\chi(Z)+(r-1)\,\chi(M).
    \end{equation*}
\end{remark}

The corresponding Fulton--Johnson identity is
\begin{equation}\label{f:FJ}
    i_*c^{FJ}(Z) = p_*c^{FJ}(Z(\hat s))-(r-1) c(TM)\cap[M].
\end{equation}
We reduce it to a universal identity.

Let $h$ denote the first Chern class of ${\mathcal O}_{\mathbb{P}(E^{\vee})}(1)$ in $H^2(\mathbb{P}(E^\vee),\Z)$. Since $\hat s$ is a regular section of $\mathcal{O}_{\mathbb{P}(E^\vee)}(1)$, $[Z(\hat s)]=h \cap [\mathbb{P}(E^\vee)]$ and, by the generalized Euler sequence \cite[B.5.8]{Ful},
\begin{equation*}
    c(T\mathbb{P}(E^\vee))=c(p^*TM) c(p^*E^{\vee}\otimes{\mathcal O}_{\mathbb{P}(E^\vee)}(1)).
\end{equation*}
The projection formula gives
\begin{equation*}
    p_*c^{FJ}(Z(\hat s)) = c(TM) \cap p_* \left(c(p^*E^{\vee}\otimes{\mathcal O}_{\mathbb{P}(E^\vee)}(1)) \frac{h}{1+h}\cap[\mathbb{P}(E^\vee)]\right).
\end{equation*}

Since $i_*[Z]=c_r(E)\cap[M]$, we also have $i_*c^{FJ}(Z)=c(TM)c(E)^{-1}c_r(E)\cap[M]$. Hence \eqref{f:FJ} is equivalent to the identity
\begin{equation}\label{f:universal}
    p_* \left(c(p^*E^{\vee}\otimes{\mathcal O}_{\mathbb{P}(E^\vee)}(1)) \frac{h}{1+h}\cap[\mathbb{P}(E^\vee)]\right) = \left(r-1+\frac{c_r(E)}{c(E)}\right)\cap[M].
\end{equation}

Aluffi proves an equivalent identity by direct computation
\cite[Lemma~3.2]{AluffiCSM}. We use universality, in a form that
will also apply to Hirzebruch classes.

\begin{lemma}\label{l:D_zero}
    Let $R$ be a commutative ring, $r\geq1$, and let $D\in R[[c_1,\dots,c_r]]$, with $\deg c_i=2i$. For $N>r$ let $G_N=\operatorname{Gr}_\C(N-r,N)$ and let $Q_N$ be its universal quotient bundle, of rank $r$. If $D(c_1(Q_N),\dots,c_r(Q_N))=0$ in $H^*(G_N;R)$ for every $N>r$, then $D=0$.
\end{lemma}
\begin{proof}
    Let $S_N$ be the tautological subbundle, so that $c(S_N)\,c(Q_N)=1$. The cohomology ring $H^*(G_N;\Z)$ is a free $\Z$-module, generated as a ring by $c_1(Q_N),\dots,c_r(Q_N)$, and
    \[
    H^*(G_N;\Z)\cong\Z[c_1,\dots,c_r]\big/\bigl(\{c(Q)^{-1}\}_{N-r+1},\dots,\{c(Q)^{-1}\}_{N}\bigr),
    \]
    where $\{c(Q)^{-1}\}_j$ denotes the component of degree $2j$ of $(1+c_1+\dots+c_r)^{-1}$; the relations express $c_j(S_N)=0$ for $j>N-r$ (see, e.g., \cite[Example~14.6.6]{Ful}). Since $H^*(G_N;\Z)$ is free, the same presentation holds with coefficients in $R$. All the relations have degree at least $2(N-r+1)$. Hence, for every $k$, the map $R[c_1,\dots,c_r]_{2k}\to H^{2k}(G_N;R)$ is injective as soon as $N-r\geq k$. Applying this to the homogeneous component of $D$ of degree $2k$, with $N$ large, shows that it vanishes. Thus $D=0$.
\end{proof}

\begin{lemma}\label{l:universal}
    The identity \eqref{f:universal} holds for all holomorphic vector bundles $E$ of rank $r$ over all compact complex manifolds $M$.
\end{lemma}
\begin{proof}
    Let $\alpha\in H_*(M,\Z)$, let $p^!$ denote the smooth homological pullback of degree $2(r-1)$, and let $a\geq0$. The projective bundle formula \cite[Proposition~3.1]{Ful} gives
    \begin{equation*}
        p_*(h^a \cap p^!\alpha) = \begin{cases}
            0, &\text{if } a < r-1 \\
            s_{a-r+1}(E^\vee) \cap \alpha, &\text{otherwise,}
        \end{cases}
    \end{equation*}
    where $s(E^\vee)=c(E^\vee)^{-1}$. Thus, for each $r$, the difference in \eqref{f:universal} is given by a universal series $D_r\in\mathbb Z[[c_1,\dots,c_r]]$, with $\deg c_i=2i$:
    \begin{align*}
        p_* \left(c(p^*E^{\vee}\otimes{\mathcal O}_{\mathbb{P}(E^\vee)}(1)) \frac{h}{1+h}\cap[\mathbb{P}(E^\vee)]\right) - \left(r-1+\frac{c_r(E)}{c(E)}\right)\cap[M]  \\
        = \text{trunc}_{2\dim M}(D_r)(c_1(E), \dots, c_r(E)) \cap [M],
    \end{align*}
    where $\text{trunc}_{n}$ is the truncation morphism from power series to polynomials of degree at most $n$.

    Suppose that $M$ is projective and $E$ is globally generated. By
    the Bertini--Kleiman theorem \cite{Kleiman}, a general section $s$ is transverse
    to the zero section. Both $Z(s)$ and $Z(\hat s)$ are then
    nonsingular by Lemma~\ref{l:fibers}(3), so their CSM and
    Fulton--Johnson classes agree. Equation \eqref{f:SM} gives \eqref{f:FJ}, hence
    $c(TM)D_r(c(E))\cap[M]=0$. Since $c(TM)$ is invertible,
    $D_r(c(E))=0$ in $H^*(M)$, with the truncation understood.

    Apply this to the projective manifold $M=G_N$ and its universal
    quotient bundle $Q_N$.
    The tautological surjection $G_N\times\mathbb C^N\to Q_N$
    shows that $Q_N$ is globally generated. Hence
    $D_r(c_1(Q_N),\dots,c_r(Q_N))=0$ for every $N>r$.
    Lemma~\ref{l:D_zero} gives $D_r=0$.
\end{proof}

\begin{theorem}{\cite[Corollary~3.6]{AluffiCSM}}\label{t:final}
    Let $s$ be a regular section of the holomorphic vector bundle $E$ of rank $r$ over the compact complex manifold $M$, and let $Z(\hat s)\subset\mathbb{P}(E^{\vee})$ be the associated hypersurface. Then, in $H_{2*}(M,\Z)$,
    \begin{equation}\label{the formula}
        i_*{\mathcal M}(Z) = p_*{\mathcal M}(Z(\hat s)),
    \end{equation}
    where the class of $Z(\hat s)$ is pushed to $\mathbb P(E^\vee)$.
\end{theorem}

\begin{proof}
    Subtract \eqref{f:SM} from \eqref{f:FJ} and use
$\dim Z(\hat s)=n+r-2\equiv n-r\pmod2$.
\end{proof}

\begin{remark}
    For $r=1$ one has $\mathbb{P}(E^{\vee}) = M$ and $\hat s=s$, and the identity \eqref{the formula} is tautological.
\end{remark}

\begin{remark}\label{r:transfer}
    Identity \eqref{the formula} transfers hypersurface formulas to
complete intersections. We apply Parusi\'nski--Pragacz in
Section~\ref{sec:closed}; applying Aluffi's hypersurface formula
\cite{Aluffi} instead gives an expression through the Segre class of
the singular scheme of $Z(\hat s)$
\cite[Corollary~2.6]{AluffiCSM}.
\end{remark}

\section{A closed formula for the Milnor class}\label{sec:closed}

Let $X$ be a hypersurface in a compact complex manifold $W$ of
dimension $m$, defined by a holomorphic section of a line bundle $L$.
For $x\in X$ let $F_x$ be the Milnor fibre at $x$ of a local equation
of $X$, and put
$$
\mu_X(x)=(-1)^{m-1}\bigl(\chi(F_x)-1\bigr),
$$
a constructible function on $X$, supported on ${\rm Sing}\,X$. With
the convention ${\mathcal M}(X)=(-1)^{\dim X}(c^{FJ}(X)-c^{SM}(X))$,
Parusi\'nski and Pragacz \cite[Theorem~0.2]{Par-Pr} proved
Yokura's conjectural formula
\begin{equationth}\label{f:PP}
{\mathcal M}(X)=c(L|_X)^{-1}\cap c_*(\mu_X).
\end{equationth}
In \cite{Par-Pr} this is written as
$\sum_S\alpha(S)\,c(L|_X)^{-1}\cap c_*(\overline S)$ over the strata of
a stratification on which $\mu_X$ is constant, with coefficients
$\alpha(S)$ defined by descending induction; since
$\sum_S\alpha(S)\mathbf{1}_{\overline S}=\mu_X$, the two forms agree.
The formula holds in the complex analytic setting.

We apply \eqref{f:PP} to $X=Z(\hat s)$, $W=\mathbb P(E^\vee)$,
$m=n+r-1$ and $L={\mathcal O}_{\mathbb P(E^\vee)}(1)$. In a local
trivialization of $E$ with $s=(f_1,\dots,f_r)$, a local equation of
$Z(\hat s)$ near $(x,[\varphi])$ is the \emph{Cayley function}
$$
g(y,\varphi)=\sum_{i=1}^r\varphi_i f_i(y),
$$
in a chart of $\mathbb P^{r-1}$. Write $\hat\mu=\mu_{Z(\hat s)}$; it is
supported on ${\rm Sing}\,Z(\hat s)$, which lies over ${\rm Sing}\,Z(s)$
by Lemma~\ref{l:fibers}(3), and $(-1)^{m-1}=(-1)^{n-r}$.

\begin{theorem}\label{t:closed}
Let $s$ be a regular section of a rank-$r$ holomorphic vector bundle
$E$ over the compact complex manifold $M$, and let $Z=Z(s)$. Then
$$
i_*\,{\mathcal M}(Z)=p_*\Bigl(c\bigl({\mathcal O}_{\mathbb P(E^\vee)}(1)
\bigr)^{-1}\cap c_*(\hat\mu)\Bigr)
\qquad\text{in }H_{2*}(M,\Z),
$$
where $\hat\mu(x,[\varphi])=(-1)^{n-r}\bigl(\chi(F_{x,\varphi})-1\bigr)$
and $F_{x,\varphi}$ is the Milnor fibre of the Cayley function $g$ at
$(x,[\varphi])$.
\end{theorem}

\begin{proof}
Apply \eqref{f:PP} to $Z(\hat s)$ and push forward by $p$, using
\eqref{the formula}.
\end{proof}

The formula is local over ${\rm Sing}\,Z$: the class $c_*(\hat\mu)$
lives on ${\rm Sing}\,Z(\hat s)$, which is fibred over ${\rm Sing}\,Z$
with fibre $\mathbb P(K_x)$ over $x$, where
$K_x=\ker((ds_x)^*:E_x^\vee\to T_x^*M)$. The differential is intrinsic along $Z$, and its corank is
$r-\operatorname{rank}(ds_x)$. The factor $c(\mathcal O(1))^{-1}$
must be applied before push-forward; it cannot in general be
replaced by a class on $M$.

\subsection{Arbitrary zero schemes}\label{ss:arbitrary}
The same argument computes the Chern--Schwartz--MacPherson class of the
zero scheme of an arbitrary section, with no regularity assumption.
This is the counterpart, in terms of Milnor fibres, of Aluffi's
expression through Segre classes \cite[Theorem~2.5]{AluffiCSM}.

\begin{theorem}\label{t:arbitrary}
Let $s$ be a nonzero holomorphic section of a rank-$r$ bundle $E$
over the compact connected complex manifold $M$, and let $Z=Z(s)$,
with no restriction on its codimension. Then
$$
i_*\,c^{SM}(Z)=c(TM)\,c(E)^{-1}c_r(E)\cap[M]
-p_*\Bigl(c\bigl({\mathcal O}_{\mathbb P(E^\vee)}(1)\bigr)^{-1}\cap
c_*(\hat\nu)\Bigr),
$$
where $\hat\nu(x,[\varphi])=\chi(F_{x,\varphi})-1$ on $Z(\hat s)$, and
$F_{x,\varphi}$ is the Milnor fibre of the Cayley function at
$(x,[\varphi])$.
\end{theorem}

\begin{proof}
Lemma~\ref{l:fibers}(1) does not use the regularity of $s$, so, as in the proof of Lemma~\ref{l:SM}, $p_*\,c^{SM}(Z(\hat s))=i_*\,c^{SM}(Z)+(r-1)\,c(TM)\cap[M]$. Since $Z(\hat s)$ is a hypersurface and $\hat{s}$ is a nonzero section of $\mathcal{O}_{\mathbb{P}(E^\vee)}(1)$, $[Z(\hat s)]=h\cap[\mathbb P(E^\vee)]$. The computation preceding \eqref{f:universal} and Lemma~\ref{l:universal}, which do not involve $s$, gives $p_*\,c^{FJ}(Z(\hat s))=c(TM)\,c(E)^{-1}c_r(E)\cap[M]+(r-1)\,c(TM)\cap[M]$. Finally, \eqref{f:PP} for $Z(\hat s)$ reads $c^{SM}(Z(\hat s))=c^{FJ}(Z(\hat s))-c({\mathcal O}(1))^{-1}\cap c_*(\hat\nu)$. Push forward by $p$ and subtract.
\end{proof}

If $s$ is regular, the first term is $i_*c^{FJ}(Z)$ and
Theorem~\ref{t:arbitrary} is equivalent to Theorem~\ref{t:closed}.
Every closed subscheme of $\mathbb P^n$ is the zero scheme of a section
$(F_1,\dots,F_r)$ of $E={\mathcal O}(d)^{\oplus r}$ for $d\gg0$, and
then $\mathbb P(E^\vee)=\mathbb P^n\times\mathbb P^{r-1}$. Thus the
Chern--Schwartz--MacPherson class of an arbitrary projective scheme is
computed from local Milnor fibres of the bihomogeneous equation
$\sum_i\varphi_iF_i$.

\begin{example}[A non-regular section]\label{ex:nonregular}
Let $M=\mathbb P^2$, $E={\mathcal O}(1)^{\oplus2}$ and $s=(x,x)$, so
that $Z$ is the line $\{x=0\}$, of codimension one, and $s$ is not
regular. Here $\mathbb P(E^\vee)=\mathbb P^2\times\mathbb P^1$,
${\mathcal O}(1)={\mathcal O}(1,1)$, and $\hat s=(\varphi_1+\varphi_2)\,x$,
so $Z(\hat s)$ is the union of $\{x=0\}\times\mathbb P^1$ and
$\mathbb P^2\times\{[1:-1]\}$. They meet transversally along the curve
$C=Z\times\{[1:-1]\}$, where the Milnor fibre is that of $xu$, with
$\chi=0$; thus $\hat\nu=-\mathbf{1}_C$. The first term of
Theorem~\ref{t:arbitrary} is $(1+H)^3(1+H)^{-2}H^2\cap[\mathbb P^2]=[pt]$.
On $C$, $c_1({\mathcal O}(1,1))=H$, so
$$
-p_*\Bigl((1-H)\cap\bigl(-[C]-2[pt]\bigr)\Bigr)=[Z]+[pt].
$$
Thus $i_*c^{SM}(Z)=[Z]+2[pt]$, the Chern class of the line.
\end{example}

\subsection{L\^e classes}
For a hypersurface $X$ in a projective manifold $W\subset\mathbb P^N$,
defined by a section of an arbitrary line bundle $L$, put
$\beta=\mu_X=(-1)^{\dim X}(\chi(F_x)-1)$. Following
\cite[Definition~1.3]{BMS-err}, define
\[
\Lambda_i(X)=\operatorname{MP}^{\mathbb P}_i(\beta)
=\sum_S(-1)^{\dim S}\eta(S,\beta)\,[P_i(\overline S)],
\]
where $\eta(S,\beta)$ is the normal Morse index and $P_i(\overline S)$
is the $i$-dimensional projective polar variety for a general flag.
These classes depend on the embedding. For reduced $X$ they are the
global L\^e classes of \cite{BMS-err}, with local representatives
given by Massey's L\^e cycles \cite{Massey}. With $H$ the hyperplane
class one has
\begin{equationth}\label{f:BMS}
{\mathcal M}_k(X)=\sum_{j\geq0}\sum_{i\geq k+j}(-1)^{i+j}
\binom{i+1}{k+j+1}\,c_1(L|_X)^j\,H^{\,i-k-j}\cap\Lambda_i(X)
\end{equationth}
\cite[Theorem~1.1]{BMS-err}, which corrects \cite[Theorem~1]{BMS}. Here $\mathcal M_k$ denotes the component of dimension $k$, with our
sign convention, also used in \cite[Example~1.8]{BMS-err}. Although
\cite[Theorem~1.1]{BMS-err} is stated for reduced sections, the same
formula holds with the above definition for nonreduced hypersurfaces:
apply \cite[Theorem~1.4]{BMS-err} to the constructible function
$\beta$ and multiply by $c(L|_X)^{-1}$, as in \eqref{f:PP}.

If $M$ is projective, then so is $\mathbb P(E^\vee)$: for an ample line
bundle $A$ on $M$, the bundle ${\mathcal O}_{\mathbb P(E^\vee)}(1)\otimes
p^*A^{m}$ is very ample for $m\gg0$. Since $p_*$ preserves dimensions,
Theorem~\ref{t:final} and \eqref{f:BMS} give the following.

\begin{corollary}\label{c:Le}
Let $M$ be projective, fix an embedding $\mathbb P(E^\vee)\subset\mathbb
P^N$ with hyperplane class $H$, and put
$h=c_1({\mathcal O}_{\mathbb P(E^\vee)}(1))$. Let $\Lambda_i$ be the
global L\^e classes of $Z(\hat s)$ defined above from the Milnor
function of its local Cayley equations $g=\sum_i\varphi_if_i$. Then, for every
$k$,
$$
i_*\,{\mathcal M}_k(Z)=p_*\Bigl(\sum_{j\geq0}\sum_{i\geq k+j}(-1)^{i+j}
\binom{i+1}{k+j+1}\,h^j\,H^{\,i-k-j}\cap\Lambda_i\Bigr).
$$
\end{corollary}

This realizes the reduction to hypersurface L\^e classes proposed
in \cite{CBMS}.

\begin{example}[Two planes meeting along a line]\label{ex:planes}
Let $M=\mathbb P^4$, $E={\mathcal O}(2)\oplus{\mathcal O}(1)$ and
$s=(x_0x_1,\ell)$, with $\ell$ a linear form independent of
$x_0,x_1$. Then $Z=Z(s)$ is the union of two planes in the hyperplane
$\{\ell=0\}$, meeting along the line $L=\{x_0=x_1=\ell=0\}$, and
${\rm Sing}\,Z=L$. Let $H$ be the hyperplane class.

\emph{Direct computation.} By inclusion--exclusion,
$c^{SM}(Z)=2\,c^{SM}(\mathbb P^2)-c^{SM}(\mathbb P^1)
=2[\mathbb P^2]+5[\mathbb P^1]+4[pt]$. Since $[Z]=2H^2\cap[\mathbb P^4]$
and $c(T\mathbb P^4)c(E)^{-1}=(1+H)^4(1+2H)^{-1}=1+2H+2H^2+\cdots$,
$c^{FJ}(Z)=2[\mathbb P^2]+4[\mathbb P^1]+4[pt]$. Hence
${\mathcal M}(Z)=c^{FJ}(Z)-c^{SM}(Z)=-[L]$.

\emph{Through Theorem~\ref{t:closed}.} At $x\in L$, $d(x_0x_1)$
vanishes and $d\ell\neq0$, so $K_x$ is the line of forms vanishing on
the second summand, and ${\rm Sing}\,Z(\hat s)=\hat L$ is a section of
$p$ over $L$. In the chart $\varphi_1=1$, $u=\varphi_2$, the Cayley
function is $g=x_0x_1+u\ell$, a nondegenerate quadratic form in four of
the five local coordinates of $\mathbb P(E^\vee)$. Its Milnor fibre at
a point of $\hat L$ is homotopic to $S^3$, so $\chi=0$ and
$\hat\mu=(-1)^{2}(0-1)=-1$ on $\hat L$. Thus
$c_*(\hat\mu)=-c(T\hat L)\cap[\hat L]=-([\hat L]+2[pt])$. Along $\hat L$
the tautological line is the summand ${\mathcal O}(-2)$ of $E^\vee$, so
${\mathcal O}(1)|_{\hat L}={\mathcal O}(2)|_L$ and
$c({\mathcal O}(1))^{-1}=1-2H$ there. Hence
$$
p_*\Bigl(c({\mathcal O}(1))^{-1}\cap c_*(\hat\mu)\Bigr)
=-(1-2H)\cap\bigl([L]+2[pt]\bigr)=-[L],
$$
in agreement with the direct computation.

\emph{Through Corollary~\ref{c:Le}.} The singular set $\hat L$ is a
single stratum, with transversal type an $A_1$ singularity of a
threefold; its normal Morse index for the Milnor function is $-1$.
Computing the MacPherson cycles of \cite{BMS-err} gives
$\Lambda_1=[\hat L]$ and $\Lambda_0=(2e-2)[pt]$, where $e$ is the
degree of $\hat L$ in $\mathbb P^N$ and $2e-2$ the degree of its polar
locus. Corollary~\ref{c:Le} gives ${\mathcal M}_1=-[L]$, and in
dimension zero
$$
\Lambda_0-2H\cdot\Lambda_1+h\cdot\Lambda_1=(2e-2)-2e+2=0,
$$
since $h\cdot\hat L=2$. The dependence on the embedding cancels in both components.
\end{example}

\section{Isolated singularities}\label{sec:examples}

\subsection{The singular locus over an isolated singular point}

Suppose that $Z$ has isolated singularities $x_1,\dots,x_N$. Then
$\mathcal M(Z)=\sum_j\mu(Z,x_j)[x_j]$ \cite{SS,Suwa}, and
Lemma~\ref{l:fibers}(3) gives
\begin{equation*}
    {\rm Sing}\,Z(\hat s)=\bigsqcup_{j=1}^N \left\{\varphi \in p^{-1}(x_j) \bigg| \sum_{i} \varphi_i df_i(x_j) = 0 \right\}.
\end{equation*}
Over $x_j$ this is the projective space $S_{x_j}\simeq\mathbb P^{c_j-1}$,
where $c_j=r-\operatorname{rank}(ds_{x_j})$. Taking degrees in
Theorem~\ref{t:final} gives 
\begin{equationth}\label{f:isolated}
    \sum_{j=1}^N\mu(Z,x_j)=\deg{\mathcal M}_0(Z(\hat s)).
\end{equationth}
To identify the contribution over each point separately, we use a
local calculation; the ambient push-forward need not be injective.

\subsection{Corank one}

In corank one, $S_x$ is a point. The equality of Milnor numbers
follows directly from the local equation.

\begin{proposition}\label{p:corank-one}
    Let $x$ be an isolated singular point of $Z=Z(s)$ of corank one.
Then $Z(\hat s)$ has a unique singular point $\xi$ over $x$, which
is isolated, and
    \begin{equation*}
        \mu(Z(\hat s),\xi)=\mu(Z,x).
    \end{equation*}
\end{proposition}

\begin{proof}
Let $s=(f_1,\dots,f_r)$ near $x$. After a change of frame of $E$ we may assume that $df_1(x),\dots,df_{r-1}(x)$ are independent and that $df_r(x)=0$. Choose coordinates $(y_1,\dots,y_{r-1},w_1,\dots,w_{n-r+1})$ centred at $x$ with $y_i=f_i$. Then $Z=\{y=0,\ g(w)=0\}$ with $g(w)=f_r(0,w)$, a hypersurface in the nonsingular germ $\{y=0\}$, so that $\mu(Z,x)=\mu(g)$.

The unique point over $x$ is $\xi=(x,[0:\dots:0:1])$. In the chart $\varphi_r=1$, with $u_i=\varphi_i$, the hypersurface $Z(\hat{s})$ is defined by
\begin{equation*}
    G=f_r(y,w)+\sum_{i<r}u_iy_i .
\end{equation*}
Write $f_r(y,w)=g(w)+\sum_{i<r}y_i\,a_i(y,w)$, where $a_i(0,0)=\partial f_r/\partial y_i(x)=0$. Then $v_i=u_i+a_i(y,w)$, together with $y$ and $w$, are coordinates at $\xi$, in which
\begin{equation*}
    G=g(w)+\sum_{i<r}y_iv_i .
\end{equation*}
This is the sum of $g$ and a nondegenerate quadratic form in separate variables. Hence $\xi$ is an isolated critical point of $G$ and $\mu(G)=\mu(g)$.
\end{proof}

\subsection{A local theorem}

\begin{definition}
    Let $x$ be an isolated singular point of $Z$. Choose local
    coordinates on a neighbourhood $V\subset M$ and a trivialization
    $s|_V=f=(f_1,\dots,f_r):V\to\mathbb C^r$.
    A sufficiently small closed \emph{Milnor ball} $B\subset V$
    centred at $x$ has the following properties \cite{Hamm}:
    \begin{itemize}
    \item $x$ is the only singular point of $Z\cap B$, and
    $Z\setminus\{x\}$ is transverse to $\partial B$;
    \item $Z\cap B$ is contractible;
    \item for every sufficiently small regular value $t$, the fibre
    $F_t=f^{-1}(t)\cap B$ is a compact manifold with boundary,
    transverse to $\partial B$, and has the homotopy type of a
    bouquet of $\mu(Z,x)$ spheres of dimension $n-r$.
    \end{itemize}
\end{definition}

Set $U=p^{-1}(B)\simeq B\times\mathbb P^{r-1}$. The section
$s_t=f-t$ has zero scheme $F_t$ and associated hypersurface
\[
X_t=Z(\hat s_t)
=\{(y,[\varphi])\in U:\varphi(f(y)-t)=0\}.
\]
Write $q_t:X_t\to B$ for the projection. Thus
$X_0=Z(\hat s)\cap U$ and $F_0=Z\cap B$.

\begin{theorem}\label{t:local}
Let $x$ be an isolated singular point of $Z$. With the notation above:
\begin{enumerate}
    \item ${\rm Sing}\,X_0=S_x$ and $\chi(X_0)=r$;
    \item for every sufficiently small regular value $t$ of $f|_B$, the hypersurface $X_t$ is nonsingular and $\chi(X_t)=r+(-1)^{n-r}\mu(Z,x)$.
\end{enumerate}
Consequently
\begin{equation*}
    \chi(X_t)-\chi(X_0)=(-1)^{n-r}\,\mu(Z,x).
\end{equation*}
\end{theorem}

\begin{proof}
Lemma~\ref{l:fibers}(3) gives $\operatorname{Sing}X_0=S_x$ and
shows that $X_t$ is nonsingular for regular $t$.
Use compactly supported Euler characteristic on locally closed
pieces; it agrees with ordinary Euler characteristic on the compact
spaces $B$, $F_t$ and $X_t$. By Lemma~\ref{l:fibers}(1),
\begin{align*}
\chi(X_t)
&=r\chi(F_t)+(r-1)\chi_c(B\setminus F_t)\\
&=\chi(F_t)+(r-1)\chi(B).
\end{align*}
Now $\chi(B)=\chi(F_0)=1$ and
$\chi(F_t)=1+(-1)^{n-r}\mu(Z,x)$ for regular $t$.
\end{proof}

Let $c=r-\operatorname{rank}(ds_x)$, let
$S_x=\mathbb P(K_x)\simeq\mathbb P^{c-1}$, and write
$H=c_1(\mathcal O_{S_x}(1))$. We use the Parusi\'nski--Pragacz
representative of the hypersurface Milnor class on $S_x$:
\begin{equation}\label{f:local-class}
\mathfrak M_x=(1+H)^{-1}\cap c_*^{S_x}(\hat\mu|_{S_x})
\quad\in H_{2*}(S_x,\mathbb Z).
\end{equation}
Here $c_*^{S_x}$ is the Chern--Schwartz--MacPherson transformation
on $S_x$. This definition uses the germ of the Cayley hypersurface
along $S_x$ and the restriction of its tautological line bundle.

\begin{lemma}\label{l:local-degree}
The degree of $(\mathfrak M_x)_0$ is $\mu(Z,x)$ and agrees with
Parusi\'nski's generalized Milnor number $\mu(Z(\hat s),S_x)$.
\end{lemma}

\begin{proof}
For a positive-dimensional ICIS, finite contact determinacy
\cite{ContactDeterminacy} gives formal equivalence of $f$ and any
system with the same sufficiently high jet. Artin approximation
\cite{Artin} makes this equivalence convergent: apply it to the
equations expressing one system as the other composed with a
coordinate change and multiplied by a matrix, preserving the
invertible linear part and the invertible constant matrix.
In dimension zero, if $I=(f_1,\ldots,f_r)$, then
$\mathfrak m^{a}\subset\mathfrak m I$ for $a$ large, so a perturbation
in $(\mathfrak m^{a})^r$ is multiplication by an invertible matrix
congruent to the identity modulo $\mathfrak m$.

Fix such a jet, of order $N$, at a point $o\in\mathbb P^n$ in an
affine chart. For $D>N$, take $r$ homogeneous polynomials of degree
$D$ with that jet in the standard affine trivialization. The space
of variations in each equation generates $\mathcal O(D)$ away
from $o$: at any point other than $o$ one can choose a linear form
vanishing at $o$ and nonzero at the given point and use its $D$th
power. The universal evaluation on the complement of $o$ is
therefore transverse to the zero section. Generic smoothness gives a
choice for which the projective zero scheme $Z'$ is smooth of the
expected codimension away from $o$. At $o$ its germ is analytically
equivalent to $(Z,x)$ by the choice of jet. In particular the
defining section is regular everywhere.

The change of coordinates and generators induces a biholomorphism of
neighbourhoods of the singular fibres in the two projective bundles.
For a change of generators by a matrix $A$, the covectors transform
by $A^{-\vee}$, preserving evaluation and identifying the tautological
lines. It follows that the Milnor functions and the classes $H$ agree
under this identification. Thus the class on the singular fibre of
the projective Cayley hypersurface $X'$ is precisely $\mathfrak M_x$.

The singular set of $X'$ is this single fibre by
Lemma~\ref{l:fibers}(3). By \eqref{f:PP} its class pushes forward to
$\mathcal M(X')$. Taking degrees in Theorem~\ref{t:final} gives
$\deg(\mathfrak M_x)_0=\deg\mathcal M_0(Z')=\mu(Z,x)$.
The degree of $\mathcal M_0(X')$ is also the sum of
Parusi\'nski's generalized Milnor numbers
\cite[Proposition~1.6]{Par}, here a single term. That number is
defined locally along the singular component, so the identification
of germs gives $\mu(Z(\hat s),S_x)=\mu(Z,x)$.
\end{proof}

\begin{theorem}[The local class and generic systems]\label{c:local}
Let $x$ be an isolated singular point of $Z(s)$, of corank $c$.
Choose a sufficiently general flag of subspaces
\[
L_0=(\mathbb C^r)^\vee\supset L_1\supset\cdots\supset L_{c-1},
\qquad\operatorname{codim}L_k=k,
\]
in a local trivialization of $E$. Restricting the linear forms in
$L_k$ to $s$ gives a system $f^{(k)}$ of $r-k$ independent linear
combinations of $f_1,\dots,f_r$. Its zero scheme has an isolated
complete intersection singularity at $x$, of dimension $n-r+k$.
Put
\[
\mu_k=\mu(Z(f^{(k)}),x),\qquad 0\leq k\leq c-1;
\quad\text{in particular, }\mu_0=\mu(Z(s),x).
\]
Then
\begin{equation}\label{f:local-full-class}
\mathfrak M_x=\mu_0[pt]+
\sum_{k=1}^{c-1}(-1)^k\mu_k\,
(1+H)^{k-1}\cap[\mathbb P^k]
\quad\text{in }H_{2*}(S_x,\mathbb Z),
\end{equation}
where $[\mathbb P^k]$ denotes the class of a linear subspace of $S_x$.
Equivalently, its characteristic numbers satisfy
\begin{equation}\label{f:local-moments}
\deg\left(\frac{H^k}{(1+H)^k}\cap\mathfrak M_x\right)
=(-1)^k\mu_k,\qquad 0\leq k\leq c-1.
\end{equation}
In particular,
\[
\mu\bigl(Z(\hat s),S_x\bigr)
=\deg(\mathfrak M_x)_0=\mu(Z(s),x),
\qquad
(\mathfrak M_x)_{c-1}=(-1)^{c-1}\mu_{c-1}[S_x].
\]
\end{theorem}

\begin{proof}
The numerical equality is Lemma~\ref{l:local-degree}. Choose a
Whitney stratification of the Cayley hypersurface near the compact
fibre over $x$ such that $S_x$ is a union of strata. The stratification
satisfies Thom's $a_g$ condition for the local defining functions $g$
of the hypersurface \cite{BrianconMaisonobeMerle}. In particular the
Milnor function is constructible for the induced stratification of
$S_x$.

A slice $B\times\mathbb P(L)$ is transverse to a stratum $A$
exactly when $\mathbb P(L)$ is transverse to the map
$A\to\mathbb P^{r-1}$ induced by projection. The action of
$GL_r(\mathbb C)$ on $\mathbb P^{r-1}$ is transitive, so parametric
transversality \cite{Kleiman} gives a general translate of a flag transverse to each
of the finitely many strata meeting a compact neighbourhood of the
fibre. The slices are then also transverse to the strata of $S_x$
inside $S_x$. For such a flag,
\[
S_k=S_x\cap\mathbb P(L_k)\simeq\mathbb P^{c-1-k}.
\]
The restricted Cayley hypersurface is the hypersurface associated to
$f^{(k)}$. It has no singular points away from $S_k$: outside $S_x$
the original hypersurface is smooth and the flag is transverse to it.
By the differential calculation in Lemma~\ref{l:fibers}(3), the zero
scheme of $f^{(k)}$ is therefore smooth away from $x$. Moreover, a
general $L_k$ meets $K_x$ in dimension $c-k$, so the differential of
$f^{(k)}$ at $x$ has rank $r-c$. For $k>0$, every component of its zero
scheme has dimension at least $n-r+k>0$, hence meets the punctured
neighbourhood, where its codimension is $r-k$. The equations thus
have the expected codimension and form a regular sequence in the
regular local ambient ring. For $k=0$ regularity is an assumption.
Each $f^{(k)}$ therefore defines an isolated complete intersection
singularity at $x$.

Let $\iota_k:S_k\hookrightarrow S_x$ be the inclusion, and let
$\hat\mu_k$ be the Milnor function of the restricted Cayley
hypersurface. Transversality to the $a_g$ stratification gives
non-characteristic restriction of the local vanishing cycles and
identifies the Milnor cohomology before and after taking the slice. The dimension changes by $k$, so
\[
\hat\mu_k=(-1)^k\iota_k^*(\hat\mu|_{S_x}).
\]
The normal bundle of $S_k$ in $S_x$ is
$\mathcal O_{S_k}(1)^{\oplus k}$. The non-characteristic
Verdier--Riemann--Roch formula for Chern classes
\cite[Theorem~1.4]{SchTrans} therefore gives
\[
c_*^{S_k}(\hat\mu_k)
=(-1)^k(1+H)^{-k}\cap
\iota_k^!c_*^{S_x}(\hat\mu|_{S_x}).
\]
After multiplication by $(1+H)^{-1}$ and push-forward to $S_x$,
this becomes
\[
(\iota_k)_*\mathfrak M_x^{(k)}
=(-1)^k\frac{H^k}{(1+H)^k}\cap\mathfrak M_x,
\]
where $\mathfrak M_x^{(k)}$ is the class \eqref{f:local-class} for
$f^{(k)}$. Applying Lemma~\ref{l:local-degree} to $f^{(k)}$ gives
$\deg(\mathfrak M_x^{(k)})_0=\mu_k$ and proves
\eqref{f:local-moments}.

These $c$ characteristic numbers determine a class in
$H_{2*}(\mathbb P^{c-1},\mathbb Z)$, because their matrix in the basis
$[pt],[\mathbb P^1],\dots,[\mathbb P^{c-1}]$ is triangular with diagonal
entries $1$. To invert it, for $k\geq1$ put
$B_k=(1+H)^{k-1}\cap[\mathbb P^k]$. For $\ell\geq1$,
\[
\deg\left(\frac{H^\ell}{(1+H)^\ell}\cap B_k\right)
=\begin{cases}1,&k=\ell,\\0,&k\ne\ell,\end{cases}
\qquad\deg B_k=0.
\]
Indeed, when $k>\ell$ the relevant coefficient is that of
$H^{k-\ell}$ in $(1+H)^{k-1-\ell}$, a polynomial of smaller degree;
when $k<\ell$ the cap product vanishes. Thus $[pt],B_1,\dots,B_{c-1}$
is the dual basis for these characteristic numbers, yielding
\eqref{f:local-full-class}.
\end{proof}

For corank one this recovers the single class $\mu_0[pt]$. For
coranks two and three it reads, respectively,
\[
\mathfrak M_x=\mu_0[pt]-\mu_1[\mathbb P^1],\qquad
\mathfrak M_x=\mu_0[pt]+(\mu_2-\mu_1)[\mathbb P^1]
+\mu_2[\mathbb P^2].
\]
Thus the higher components recover the Milnor numbers of generic
systems with fewer equations. If $Z(s)$ has several isolated
singularities, the formula applies separately over each point.

\begin{remark}\label{r:invariance}
The moments \eqref{f:local-moments} give independence of the
general flag. A change of local generators induces an isomorphism of
the Cayley germs preserving their tautological lines, as in
Lemma~\ref{l:local-degree}.

To see independence of the embedding, put $b=r-c$. By the implicit
function theorem and a change of generators the system can be written
$(y_1,\ldots,y_b,g_1(z),\ldots,g_c(z))$, where $g$ gives the
presentation in the smooth germ $y=0$ and $dg(0)=0$. In a chart over
$S_x$ where one of the coefficients of the $g_j$ is $1$, the
Cayley equation is the Cayley equation for $g$ plus
$\sum_{i=1}^b u_i y_i$. Thom--Sebastiani preserves the signed Milnor
function under this suspension by a nondegenerate quadratic form in
$2b$ variables. The identification of the singular fibres preserves
$H$, so \eqref{f:local-class} and all its moments agree with those of
the presentation in $y=0$.

This last embedding is minimal, since its ambient dimension is the
embedding dimension of $(Z,x)$. An isomorphism between two such
germs lifts to a biholomorphism of their minimal smooth ambient
germs: lift the coordinate functions and use that they induce an
isomorphism on cotangent spaces. Two regular sequences generating
the same ideal are related by an invertible matrix over the local
ring. Since $H_{2*}(S_x,\mathbb Z)$ is generated by the classes of linear
subspaces, which are matched by these identifications, $\mathfrak M_x$
and all the $\mu_k$ depend only on the analytic germ.
\end{remark}

\begin{remark}\label{r:LeGreuel}
The flag gives nested systems: $f^{(k+1)}$ consists of $r-k-1$ of the
equations in $f^{(k)}$. For $0\leq k\leq c-2$, both are isolated
complete intersection singularities, and the L\^e--Greuel formula
\cite{Le-CI,Greuel} gives
\[
\mu_k+\mu_{k+1}
=\dim_{\mathbb C}\mathcal O_{\mathbb C^n,x}\big/
\bigl(f^{(k+1)},J_{r-k}(f^{(k)})\bigr),
\]
where $J_{r-k}$ is the ideal of maximal minors of the Jacobian matrix.
Thus \eqref{f:local-full-class} expresses $\mathfrak M_x$ by $\mu_{c-1}$
and these colengths. For $f=(x^2,y^2)$ the colength is
$\dim\mathcal O/(x^2+y^2,xy)=4=\mu_0+\mu_1$. For
$f=(x^2,y^2,z^2)$ it is the intersection multiplicity of the quartic
curve $f^{(1)}=0$ with the cubic $xyz=0$, namely $12=\mu_0+\mu_1$.
\end{remark}

\subsection{Examples}

\begin{example}\label{ex:conic}
Let $M=\mathbb{P}^2$, $E={\mathcal O}(2)\oplus{\mathcal O}(1)$ and $s=(q,\ell)$, where $q=0$ is a smooth conic and $\ell=0$ is a line tangent to it at a point $x$. Then $Z(s)$ is the point $x$ with multiplicity two, so $c^{SM}(Z)=[x]$, $c^{FJ}(Z)=2[x]$ and ${\mathcal M}(Z)=[x]$. Here $\mathbb{P}(E^{\vee})$ is a $\mathbb{P}^1$-bundle over $\mathbb{P}^2$ and $Z(\hat s)$ is a surface. In affine coordinates with $q=y-x^2$ and $\ell=y$, the surface has one singular point, at $[\varphi]=[1:-1]$ over the origin, where it is given by $uy-x^2=0$: an $A_1$ singularity. Thus ${\mathcal M}(Z(\hat s))$ is the class of that point, and $p_*{\mathcal M}(Z(\hat s))=[x]={\mathcal M}(Z)$. Also $\chi(Z(\hat{s}))=4=\chi(Z)+\chi(\mathbb{P}^2)$, in agreement with \eqref{f:SM}: the map $q:Z(\hat{s}) \to M$ is an isomorphism over $\mathbb{P}^2\setminus\{x\}$ and its fibre over $x$ is a line.
\end{example}

\begin{example}\label{ex:corank-two}
Let $f=(x^2,y^2)$ on $\C^2$, so that $n=r=2$ and $Z=f^{-1}(0)$ is the
origin with multiplicity four, a complete intersection singularity of
dimension zero and corank two, with $\mu(Z,0)=3$: the Milnor fibre
consists of four points. Here $K_0=\C^2$, and the hypersurface
$$
X_0=\{\varphi_1x^2+\varphi_2y^2=0\}\subset B\times\mathbb{P}^1
$$
is singular along the whole line $S_0=\{0\}\times\mathbb{P}^1$. It
has an ordinary double curve along $S_0$ away from the two points
$[1:0]$ and $[0:1]$, where it has pinch points. By
Theorem~\ref{t:local}, $\chi(X_0)=2$ and $\chi(X_t)=1+4=5$, so that
$\mu(X_0,S_0)=3=\mu(Z,0)$. A general linear combination of the two
equations is a nondegenerate quadratic form, with $\mu_1=1$.
Theorem~\ref{c:local} therefore gives the full local class
\[
\mathfrak M_0=-[\mathbb P^1]+3[pt].
\]
More generally, for $f=(x^a,y^b)$ with $a,b\geq2$, one has
$\mu_0=ab-1$ and $\mu_1=(a-1)(b-1)$, so
\[
\mathfrak M_0=-(a-1)(b-1)[\mathbb P^1]+(ab-1)[pt].
\]
For a corank-three example, take $f=(x^2,y^2,z^2)$ on $\mathbb C^3$.
Here $\mu_0=7$, $\mu_1=5$ and $\mu_2=1$: a generic two-equation
smoothing compactifies to a smooth intersection of two quadrics in
$\mathbb P^3$, an elliptic curve with four points removed, while a
generic single equation is a nondegenerate quadratic form.
Thus $\mathfrak M_0=[\mathbb P^2]-4[\mathbb P^1]+7[pt]$.
Independently, its Milnor function on $S_0=\mathbb P^2$ is
$1$ off the coordinate lines, $-1$ on their one-dimensional strata,
and $1$ at the three vertices. These values follow by
Thom--Sebastiani from the local blocks $z^2$ and $uz^2$.
Consequently
$\hat\mu=\mathbf{1}_{\mathbb P^2}-2\sum_{i=1}^3\mathbf{1}_{L_i}
+4\sum_{j=1}^3\mathbf{1}_{p_j}$, so
\[
(1+H)^{-1}\cap c_*(\hat\mu)
=(1+H)^{-1}\cap
\bigl([\mathbb P^2]-3[\mathbb P^1]+3[pt]\bigr)
=[\mathbb P^2]-4[\mathbb P^1]+7[pt],
\]
in agreement with \eqref{f:local-full-class}.
\end{example}

\subsection{Revisiting the former comparison formula}\label{ss:old}

Let $0\to Q^\vee\to p^*E\to{\mathcal O}_{\mathbb{P}(E^{\vee})}(1)\to 0$ be the dualized tautological sequence \eqref{eq:dual_tautological_seq}. Since the image $\hat s$ vanishes on $Z(\hat s)$, the section $p^*s$
restricts there to a section $\tilde s$ of $Q^\vee|_{Z(\hat s)}$. The formula stated in \cite[Theorem~6.4]{CBMS} is, up to sign, 
\begin{equationth}\label{f:old}
    i_*\mathcal{M}(Z) = p_*\left(c\left(p^*E^{\vee}\otimes    \mathcal{O}_{\mathbb{P}(E^{\vee})}(1)\right)^{-1}    c_1\left(\mathcal{O}_{\mathbb{P}(E^{\vee})}(1)\right)^{r-1}c(Q^\vee)^{-1}c_{top}(Q^\vee)\cap \mathcal{M}(Z(\hat s))\right).
\end{equationth}
A 2018 revision of the manuscript explicitly imposes the assumption that $Z(\hat s)$ has a Whitney stratification adapted to its singular locus, such that $\tilde s$ is transverse to the zero section along each stratum; this condition is required to apply \cite[Proposition~1.3]{PP}. This hypothesis is absent from \cite[Theorem~6.4]{CBMS}.

\begin{proposition}
    Given a holomorphic complex vector bundle $\pi:E \to M$ of rank $r \geq 2$ and a regular section $s \in H^0(M,E)$, let $Z$ denote the zero locus of $s$, and define $\hat{s}$ and $\tilde{s}$ as above. If $Z$ is singular, then for all Whitney stratifications $\{S_\alpha\}$ of $Z(\hat{s})$ adapted to its singular locus, there exists a stratum where $\tilde{s}$ is identically zero.
\end{proposition}
\begin{proof}
    We have $Z(\tilde s)=q^{-1}(Z)$, which contains
    $\operatorname{Sing}Z(\hat s)$ by Lemma~\ref{l:fibers}(3).
    Since $Z$ is singular, this is a nonempty union of strata.
    The section $\tilde s$ vanishes identically on each of them.
\end{proof}

Since $\operatorname{rank}Q^\vee=r-1>0$, an identically zero
section on a stratum cannot be transverse to the zero section.
The 2018 transversality hypothesis therefore forces $Z(\hat s)$,
and hence $Z(s)$, to be nonsingular; both Milnor classes then vanish.

\begin{remark}
    The unrestricted formula \eqref{f:old} of \cite{CBMS} is false in general, by a dimension count. The class $\mathcal M(Z(\hat s))$ is supported on ${\rm Sing}\,Z(\hat s)$, whose dimension is at most $\dim{\rm Sing}\,Z+r-1$ by Lemma~\ref{l:fibers}. In \eqref{f:old} it is capped with $c_{top}(Q^\vee)$ and with $c_1(\mathcal O(1))^{r-1}$, which lower the dimension by $2(r-1)$, while the remaining factors only add components of lower dimension. Hence, for $r\geq2$, the right-hand side of \eqref{f:old} has no component of dimension greater than $\dim{\rm Sing}\,Z-(r-1)<\dim{\rm Sing}\,Z$. In Example~\ref{ex:conic}, $\dim{\rm Sing}\,Z=0$, so the right-hand side vanishes, whereas $\mathcal M(Z)=[x]$.
\end{remark}

\section{Hirzebruch--Milnor classes}\label{sec:hirzebruch}

We now work over a nonsingular complete complex algebraic variety
$M$, with an algebraic rank-$r$ bundle $E$ and a regular section $s$.
The fibre description of Lemma~\ref{l:fibers} gives a motivic
comparison. We use even Borel--Moore homology with rational coefficients.

\subsection{Hirzebruch classes}

Write $T_{y*}$ for the normalized motivic Hirzebruch class
transformation, often denoted $\hat T_{y*}$. Its characteristic
power series is
\begin{align}
    Q_y(\alpha) = \frac{\alpha(1+y)}{1 - e^{-\alpha (1+y)}} - \alpha y \in \mathbb{Q}[y][[\alpha]], & Q_y(0) = 1
\end{align}
and define $T_{y}^*(E) = \prod_j Q_y(\alpha_j)$, where the $\alpha_j$ are the Chern roots of $E$. In particular $T_{-1}^*(E) = c(E)$ and $T^*_{0}(E)$ is the Todd class. Since its degree-zero component is $1$, the class $T_y^*(E)$ is
invertible in $H^*(M,\mathbb Q)[y]$. This normalization is used in
the virtual comparison below; the unnormalized series has constant
term $1+y$.

Let $K_0({\rm var}/X)$ be the relative Grothendieck group of complex algebraic varieties over $X$. It is a module over $K_0({\rm var}/pt)$, with $[X'\to X]\cdot[F]=[X'\times F\to X]$, and we write $\mathbb L=[\mathbb A^1]$. The normalized form of the motivic Hirzebruch class transformation
of Brasselet, Sch\"urmann and Yokura \cite{BSY} is
\begin{equation*}
    T_{y*}:K_0({\rm var}/X)\longrightarrow H_*(X)\otimes\mathbb Q[y].
\end{equation*}
It commutes with proper push-forward and exterior products, and
satisfies $T_{y*}([{\rm id}_X])=T_y^*(TX)\cap[X]$ when $X$ is nonsingular. For $X=pt$ the transformation is the compactly supported Hodge polynomial $\chi_y$; in particular $\chi_y(\mathbb L)=-y$ and $\chi_y(\mathbb P^k)=1+(-y)+\dots+(-y)^k$. Hence
\begin{equationth}\label{f:module}
    T_{y*}\bigl([X'\to X]\cdot[F]\bigr)=\chi_y(F)\;T_{y*}([X'\to X]).
\end{equationth}

The Hirzebruch class of a variety $X$ is $T_{y*}(X)=T_{y*}([{\rm id}_X])$, and for $y=-1$ it is the Schwartz--MacPherson class $c^{SM}(X)\otimes\mathbb Q$.

If $X$ is a local complete intersection in a nonsingular variety $W$, defined by a regular section of a vector bundle $N$, its virtual Hirzebruch class is
\begin{equation*}
    T^{\rm vir}_{y*}(X)=T_y^*(TW|_X)\,T_y^*(N|_X)^{-1}\cap[X].
\end{equation*}

We use the sign convention
\begin{equationth}\label{f:HM-convention}
    {\mathcal M}_{y*}(X)=T^{\rm vir}_{y*}(X)-T_{y*}(X),
\end{equationth}
with the order of subtraction used in \cite{CMSS,MSS}. Our
notation refers to the difference on $X$; the localized class
$M_y(X)$ of \cite[Theorem~1]{MSS} lies on
$\Sigma=\operatorname{Sing}X$ and represents this difference after
push-forward.

For $y=-1$ the virtual class is the Fulton--Johnson class, so that ${\mathcal M}_{-1*}(X)=(-1)^{\dim X}{\mathcal M}(X)\otimes\mathbb Q$. The class ${\mathcal M}_{y*}(X)$ is supported on the singular locus of $X$.

\subsection{The motivic class of the associated hypersurface}

\begin{lemma}\label{l:motivic}
    In $K_0({\rm var}/M)$,
    \begin{equation*}
        [Z(\hat s)\to M]=[\mathbb P^{r-2}]\cdot[{\rm id}_M]+\mathbb L^{r-1}\cdot[Z(s)\to M],
    \end{equation*}
    where $[\mathbb P^{-1}]=0$.
\end{lemma}

\begin{proof}
Over $Z$, the hypersurface $Z(\hat s)$ is the projective bundle
$\mathbb P(E^\vee|_Z)$. Over $M\setminus Z$, it is the
projectivization of the rank-$(r-1)$ annihilator bundle of $s$.
Both are Zariski locally trivial, so
\begin{equation*}
    [Z(\hat s)\to M]=[\mathbb P^{r-1}]\cdot[Z\to M] +[\mathbb P^{r-2}]\cdot\bigl([{\rm id}_M]-[Z\to M]\bigr),
\end{equation*}
and $[\mathbb P^{r-1}]-[\mathbb P^{r-2}]=\mathbb L^{r-1}$.
\end{proof}

\subsection{The comparison}

\begin{theorem}\label{t:hirzebruch}
    Let $s$ be a regular section of the algebraic vector bundle $E$ of rank $r$ over the nonsingular complete variety $M$. Then, in $H_*(M)\otimes\mathbb Q[y]$,
    \begin{align}
        p_*\,T_{y*}(Z(\hat s))&=\chi_y(\mathbb P^{r-2})\,T_{y*}(M) + (-y)^{r-1}\,i_*\,T_{y*}(Z(s)),\label{f:Ty}\\
        p_*\,T^{\rm vir}_{y*}(Z(\hat s))&=\chi_y(\mathbb P^{r-2})\,T_{y*}(M) +(-y)^{r-1}\,i_*\,T^{\rm vir}_{y*}(Z(s)),\label{f:Tvir}\\
        p_*\,{\mathcal M}_{y*}(Z(\hat s))&=(-y)^{r-1}\;i_*\,{\mathcal M}_{y*}(Z(s)).\label{f:My}
    \end{align}
\end{theorem}

\begin{proof}
Apply $T_{y*}$ to Lemma~\ref{l:motivic}, using \eqref{f:module}
and proper covariance, to obtain \eqref{f:Ty}.

For \eqref{f:Tvir}, let $P=\mathbb P(E^\vee)$ and $h=c_1({\mathcal O}_P(1))$. Since $[Z(\hat s)]=h\cap[P]$ and $T_y^*(TP)=p^*T_y^*(TM)\cdot T_y^*(T_{P/M})$, and since $i_*[Z(s)]=c_r(E)\cap[M]$, formula \eqref{f:Tvir} is equivalent to
\begin{equationth}\label{f:universal-y}
    p_*\Bigl(T_y^*(T_{P/M})\,T_y^*({\mathcal O}_P(1))^{-1}\,h\cap[P]\Bigr)=\Bigl(\chi_y(\mathbb P^{r-2})+(-y)^{r-1}\,\frac{c_r(E)}{T_y^*(E)}\Bigr)\cap[M].
\end{equationth}

As in Lemma~\ref{l:universal}, the difference in
\eqref{f:universal-y} is a universal series
$D_y(c_1(E),\dots,c_r(E))\cap[M]$ with coefficients in $\mathbb Q[y]$.
For globally generated $E$ and general $s$, both zero schemes are
nonsingular, so their virtual and motivic classes agree. Equation
\eqref{f:Ty} therefore gives
$T_y^*(TM)D_y(c(E))\cap[M]=0$. Since $T_y^*(TM)$ is a unit,
$D_y(c(E))=0$. Evaluation on the universal quotient bundles and
Lemma~\ref{l:D_zero} give $D_y=0$.

Formula \eqref{f:My} is the difference of \eqref{f:Tvir} and \eqref{f:Ty}.
\end{proof}

For $y=-1$ one has $\chi_{-1}(\mathbb P^{r-2})=r-1$ and
$(-y)^{r-1}=1$, and Theorem~\ref{t:hirzebruch} becomes
Theorem~\ref{t:final} with rational coefficients. For $y=0$ and
$r\geq2$ it says that $p_*\,{\mathcal M}_{0*}(Z(\hat s))=0$.

\begin{example}
Consider again the conic and its tangent line of
Example~\ref{ex:conic}. The Hirzebruch class of $Z(s)$ is that of the
reduced point, $T_{y*}(Z)=[x]$, while
$T^{\rm vir}_{y*}(Z)=2[x]$, because the normalized characteristic
classes have constant term $1$ and the fundamental class of the scheme
$Z$ has degree two. Thus $\mathcal M_{y*}(Z(s))=[x]$. The associated surface has an
ordinary double point $\xi$, whose vanishing cohomology is
one-dimensional of Hodge type $(1,1)$. The isolated-hypersurface
formula \cite[Corollary~2, formula~(0.4)]{MSS}, \cite{CMSS} gives
${\mathcal M}_{y*}(Z(\hat s))=-y\,[\xi]$. Hence
$p_*{\mathcal M}_{y*}(Z(\hat s))=-y\,[x]=(-y)^{r-1}{\mathcal M}_{y*}(Z(s))$,
with $r=2$.
\end{example}

\subsection{A closed formula for Hirzebruch--Milnor classes}

We also write $T_{y*}$ for the normalized Hirzebruch class
transformation on the Grothendieck group of bounded complexes of
mixed Hodge modules \cite{MSS}. Following that paper, $\mathbb Q_W^H$ denotes the constant complex
with underlying rational complex $\mathbb Q_W$, without a dimension
shift. We use the derived vanishing-cycle functor $\varphi_h$ whose
stalk cohomology is reduced Milnor cohomology; this includes the
compensating shift from the perverse-exact convention. Classes supported on $\Sigma=\operatorname{Sing}X$ are
pushed forward to $X$ in the formulas below.

For a hypersurface $X$ in a nonsingular projective variety $W$,
defined by a section $\sigma$ of a very ample line bundle $L$ (or of a
power of one), Maxim, Saito and Sch\"urmann give the explicit formula
\begin{equationth}\label{f:MSS}
{\mathcal M}_{y*}(X)=T_{y*}\Bigl((j_\Sigma)_!\,\varphi_h\,\mathbb Q^H_W\Bigr),
\qquad h=\sigma/\sigma',
\end{equationth}
where $\sigma'$ is a general section of $L$, the function $h$ is
defined on $W\setminus Z(\sigma')$, and
$j_\Sigma:{\rm Sing}\,X\setminus Z(\sigma')\hookrightarrow{\rm Sing}\,X$
\cite[Theorem~1 and Corollary~2, formula~(0.3), case $r=1$]{MSS}. At $y=-1$ this gives \eqref{f:PP} with rational coefficients after
multiplication by $(-1)^{\dim X}$. Indeed, if $D'=Z(\sigma')$ and
$\nu(x)=\chi(F_x)-1$, the constructible function of the extension
by zero is $\nu\,\mathbf 1_{\Sigma\setminus D'}$. Additivity and
non-characteristic Verdier--Riemann--Roch for the inclusion of $D'$
give
\[
c_*(\nu\,\mathbf 1_{\Sigma\setminus D'})
=c(L|_X)^{-1}\cap c_*(\nu)
\]
after push-forward to $X$.

\begin{theorem}\label{t:closed-y}
In the setting of Theorem~\ref{t:hirzebruch}, assume that
${\mathcal O}_{\mathbb P(E^\vee)}(1)$ is very ample, for instance that
$E$ is very ample. Let $s'$ be a general section of $E$, and
$h=\hat s/\hat s'$. Then
$$
(-y)^{r-1}\,i_*\,{\mathcal M}_{y*}(Z(s))
=p_*\,T_{y*}\Bigl((j_\Sigma)_!\,\varphi_h\,\mathbb Q^H_{\mathbb P(E^\vee)}\Bigr).
$$
In particular, the right-hand side is divisible by $(-y)^{r-1}$ in
$H_*(M)\otimes\mathbb Q[y]$.
\end{theorem}

\begin{proof}
Since $H^0(\mathbb P(E^\vee),{\mathcal O}(1))=H^0(M,E)$, a general
section of ${\mathcal O}(1)$ is $\hat s'$ for a general $s'$. Apply
\eqref{f:MSS} to $X=Z(\hat s)$ and use \eqref{f:My}.
\end{proof}

Locally, $h=(\sum_i\varphi_if_i)/(\sum_i\varphi_if'_i)$.
The theorem computes the Hirzebruch--Milnor class from this function,
without induction on $r$. The very ampleness assumption comes from
the hypersurface case of \cite{MSS}; it would be useful to know
whether global generation of $E$ suffices.

\begin{remark}
In the situation of Proposition~\ref{p:corank-one}, the hypersurface
$Z(\hat s)$ is locally given by $g(w)+\sum_{i<r}y_iv_i$, the sum of
the function $g$ defining $Z(s)$ and a nondegenerate quadratic form
in $2(r-1)$ further variables. By the Thom--Sebastiani theorem for
filtered $\mathcal D$-modules and mixed Hodge modules
\cite{MSS-TS}, each block $y_iv_i$ contributes a Tate twist. After forgetting
monodromy and taking $\chi_y$, the $r-1$ blocks give the factor
$(-y)^{r-1}$ in \eqref{f:My}, the local Hodge-theoretic form of
$\mathbb L^{r-1}$ in Lemma~\ref{l:motivic}. A corresponding
Hodge-theoretic local comparison in arbitrary corank would extend
this interpretation beyond the quadratic normal form.
\end{remark}

 \end{document}